
\input amstex
\documentstyle{amsppt}

\vcorrection{-0.5cm}
\hcorrection{+0.3cm}
\pageheight{20.7cm}
\pagewidth{12.7cm}
\magnification=\magstep1
\parskip=6pt

\define\ce{{\Cal E}}
\define\cf{{\Cal F}}
\define\cg{{\Cal G}}

\define\co{{\Cal O}}
\define\cs{{\Cal S}}

\define\pp{{\Bbb P}}
\define\qq{{\Bbb Q}}

\define\zz{{\Bbb Z}}
\define\ff{{\Bbb F}}

\define\pe{{\Bbb P(\Cal E)}}
\define\he{{H(\Cal E)}}
\define\pf{{\Bbb P(\Cal F)}}
\define\hf{{H(\Cal F)}}
\define\pg{{\Bbb P(\Cal G)}}
\define\hg{{H(\Cal G)}}

\define\pic{\operatorname{Pic}}
\define\rk{\operatorname{rank}}

\topmatter 

\title  A generalization of curve genus for ample vector bundles, II
\endtitle

\author {Yoshiaki FUKUMA$^{*}$ and Hironobu ISHIHARA}
\endauthor

\address Department of Mathematics,
 Faculty of Science,
 Tokyo Institute of Technology, 
 Oh-okayama, Meguro-ku, Tokyo 152,
 Japan
\endaddress
\email fukuma\@math.titech.ac.jp
\endemail

\address Department of Mathematics,
 Faculty of Science,
 Tokyo Institute of Technology, 
 Oh-okayama, Meguro-ku, Tokyo 152,
 Japan
\endaddress
\email ishihara\@math.titech.ac.jp
\endemail

\thanks
{$^{*}$Research Fellow of the Japan Society for the Promotion of Science}
\endthanks

\abstract  
  Let $X$ be a compact complex manifold of dimension $n\ge 2$ and
  $\ce$ an ample vector bundle of rank $r<n$ on $X$.
  As the continuation of Part I, we further study the properties 
  of $g(X,\ce)$ that is an invariant for pairs $(X,\ce)$ and
  is equal to curve genus when $r=n-1$.
  Main results are the classifications of $(X,\ce)$ with 
  $g(X,\ce)=2$ (resp. 3) when $\ce$ has a regular section
  (resp. $\ce$ is ample and spanned) and $1<r<n-1$.
\endabstract


 
\leftheadtext{FUKUMA AND ISHIHARA}

\endtopmatter

\document

\subhead{Introduction}
\endsubhead

The present paper is a continuation of \cite{I}.
For a pair $(X,\Cal{E})$ which consists of a compact complex
manifold $X$ of dimension $n\geq 2$
and an ample vector bundle $\Cal{E}$ of rank $r<n$ on $X$,
we defined in \cite{I} an invariant $g(X,\Cal{E})$ by the formula
$$2g(X,\Cal{E})-2:=(K_{X}+(n-r)c_{1}(\Cal{E}))c_{1}(\Cal{E})^{n-r-1}
c_{r}(\Cal{E}).$$
We note that $g(X,\Cal{E})$ is a non-negative integer,
and $g(X,\Cal{E})$ is equal to the curve genus of $(X,\Cal{E})$ 
when $r=n-1$.
As in the case of curve genus, above $(X,\Cal{E})$ with 
$g(X,\Cal{E})\leq 1$ have been classified in \cite{I}; 
moreover, it is shown that $g(X,\Cal{E})\geq q(X)$
for spanned $\Cal{E}$ and its equality condition is given in \cite{I}.
($q(X)$ is the irregularity of $X$.)

After we recall some preliminary results in Section 1,
we consider the cases $g(X,\Cal{E})=2$ and $g(X,\Cal{E})=3$
when $1<r<n-1$ in Section 2.
Corresponding results for $c_1$-sectional genus are given in
\cite{Fj2} and \cite{BiLL} respectively.
In Section 3 we consider the cases $g(X,\Cal{E})=q(X)+1$ and 
$g(X,\Cal{E})=q(X)+2$ when $1<r<n-1$.
Related results for $c_1$-sectional genus are given in \cite{R}.
In Section 4 we give another relation between $g(X,\Cal{E})$ and $q(X)$,
namely $g(X,\Cal{E})\ge 2q(X)-1$ for $1<r<n-1$.
When $r=1$, this inequality is satisfied except one case.
In Section 5 we show that $g(X,\Cal{E})\ge g(C)$ when there exists
a fibration $f:X\to C$ over a curve.
We also give its equality condition.
Finally in Appendix we give a classification of $(X,L)$
with $g(X,L)=q(X)+2$ and $n=2$ for ample and spanned line bundles
$L$ on $X$.

The authors would like to express their gratitude to Professor
Takao Fujita for his valuable comments, especially for informing
them of Lemma 2.2.5.

\subhead{1. Preliminaries}
\endsubhead

We use a notation similar to that in \cite{I}.
For example, we denote by $H(\Cal{E})$ the tautological line bundle
on $\Bbb{P}_{X}(\Cal{E})$, the projective space bundle  
associated to a vector bundle $\Cal{E}$ on a variety $X$.
We say that a vector bundle $\Cal{E}$ is spanned 
if $H(\Cal{E})$ is spanned.
A polarized manifold $(X,L)$ is said to be a scroll over a variety $W$
if $(X,L)\simeq(\pp_W(\cf),\hf)$ for some ample vector bundle $\cf$ 
on $W$.
We denote by $\ff_e$ the Hirzebruch surfaces 
$\pp_{\pp^1}(\co_{\pp^1}\oplus\co_{\pp^1}(-e))$ ($e>0$),
by $\sigma$ the minimal section, and by $f$ a fiber of the ruling
$\ff_e\to\pp^1$.
Numerical equivalence is denoted by $\equiv$.

\definition{Definition 1.1}
  Let $X$ be a compact complex manifold of dimension $n\ge 2$
  and $\ce$ an ample vector bundle of rank $r<n$ on $X$.
  We define a rational number $g(X,\ce)$ for the pair
  $(X,\ce)$ by the formula
  $$ 2g(X,\ce)-2:=(K_X+(n-r)c_1(\ce))c_1(\ce)^{n-r-1}c_r(\ce). $$
  It turns out that $g(X,\ce)$ is a non-negative integer
  (see \cite{I}).
  When $r=1$ (resp. $r=n-1$), $g(X,\ce)$ is nothing but
  the sectional genus (resp. curve genus) of $(X,\ce)$.
\enddefinition

\example{Remark 1.2}
  Let $(X,\ce)$ be as above.
  Suppose that $(X,\ce)$ satisfies the condition
  \roster
  \item"($\ast$)" there exists a section $s\in H^0(X,\ce)$ whose zero 
            locus $Z:=(s)_0$ is a smooth submanifold of $X$
            of the expected dimension $n-r$.
 \endroster
 Then we have $g(X,\ce)=g(Z,\det\ce_Z)$ (see \cite{I}).
 If $\ce$ is spanned, then $\ce$ satisfies ($\ast$) 
 by Bertini's theorem.
\endexample

The following facts are used in the subsequent sections.

\proclaim{Proposition 1.3}
 Let $X$ be an $n$-dimensional compact complex manifold
 and $\ce$ an ample vector bundle of rank $r<n$ on $X$
 with the property $(\ast)$ in {\rm(1.2)}.
 Let $\iota:Z\hookrightarrow X$ be the embedding.
 Then
 \roster
 \item $H^i(\iota):H^i(X,\zz)\to H^i(Z,\zz)$ is an isomorphism
       for $i<n-r$.
 \item $H^i(\iota)$ is injective and its cokernel is torsion free
       for $i=n-r$.
 \item $\pic(\iota):\pic(X)\to\pic(Z)$ is an isomorphism
       for $n-r\ge 3$.
 \item $\pic(\iota)$ is injective and its cokernel is torsion free
       for $n-r=2$.
 \endroster
\endproclaim

\demo{Proof} See Theorem 1.3 in \cite{LM1} and see also Theorem 1.1
             in \cite{LM2}. \qed
\enddemo

\proclaim{Proposition 1.4}
 Let $X$ be an $n$-dimensional compact complex manifold
 and $\ce$ an ample vector bundle of rank $r\ge 2$ on $X$
 with the property $(\ast)$. 

 If $Z\simeq\pp^{n-r}(n-r\ge 1)$, then $(X,\ce)$ is one of the 
 following:
 \roster
 \item"(P1)" $(\pp^n,\co_{\pp^n}(1)^{\oplus r})$;
 \item"(P2)" $(\pp^n,\co_{\pp^n}(2)\oplus\co_{\pp^n}(1)^{\oplus(n-2)})$;
 \item"(P3)" $(\qq^n,\co_{\qq^n}(1)^{\oplus(n-1)})$;
 \item"(P4)" $X\simeq\pp_{\pp^1}(\cf)$ for some vector bundle $\cf$
             of rank $n$ on $\pp^1$ and
             $\ce=\oplus_{j=1}^{n-1}(H(\cf)\mathbreak
                                           +\pi^*\co_{\pp^1}(b_j))$,
             where $\pi:X\to\pp^1$ is the bundle projection.
 \endroster

 If $Z\simeq\qq^{n-r}(n-r\ge 2)$, then $(X,\ce)$ is one of the 
 following:
 \roster
 \item"(Q1)" $(\pp^n,\co_{\pp^n}(2)\oplus\co_{\pp^n}(1)^{\oplus(r-1)})$;
 \item"(Q2)" $(\qq^n,\co_{\qq^n}(1)^{\oplus r})$;
 \item"(Q3)" $X\simeq\pp_{\pp^1}(\cf)$ and
             $\ce=\oplus_{j=1}^{n-2}(H(\cf)+\pi^*\co_{\pp^1}(b_j))$,
             where $\cf$ is the same as that in {\rm(P4)}.
 \endroster
\endproclaim

\demo{Proof} See Theorem A and Theorem B in \cite{LM1}. \qed
\enddemo

\proclaim{Proposition 1.5}
Let $X$ be a complex projective manifold of dimension $n$ and let
$\Cal{E}$ be an ample vector bundle of rank $n-2\geq 2$ on $X$
satisfying $(\ast)$.
\roster
\item If $Z$ is a geometrically ruled surface over a smooth curve $B$
      such that $Z\not= \Bbb{F}_{0}, \Bbb{F}_{1}$, then
      $X$ is a $\Bbb{P}^{n-1}$-bundle over $B$ and 
      $\Cal{E}_{F}=\Cal{O}_{\Bbb{P}^{n-1}}(1)^{\oplus(n-2)}$
      for every fiber $F$ of the bundle map $X\to B$.
\item If $Z=\Bbb{F}_{0}$, then $(X,\Cal{E})$ is either the type in 
      {\rm(1)} with $B=\pp^1$ or
      $(\Bbb{P}^{n}, \Cal{O}_{\Bbb{P}^{n}}(2)\oplus
                     \Cal{O}_{\Bbb{P}^{n}}(1)^{\oplus(n-3)})$ 
      or
      $(\Bbb{Q}^{n}, \Cal{O}_{\Bbb{Q}^{n}}(1)^{\oplus(n-2)})$.
\item If $Z=\Bbb{F}_{1}$, then $(X,\Cal{E})$ is either the type in 
      {\rm(1)} with $B=\pp^1$ or possibly
      $X\simeq \Bbb{P}_{\Bbb{P}^{2}}(\Cal{F})$ for some ample vector 
      bundle $\Cal{F}$ on $\Bbb{P}^{2}$ with $c_{1}(\Cal{F})=k(n-2)+3$ 
      for some positive integer $k$ and 
      $\Cal{E}_{F}=\Cal{O}_{\Bbb{P}^{n-2}}(1)^{\oplus(n-2)}$
      for every fiber $F$ of the bundle map $X\to \Bbb{P}^{2}$.
\endroster
\endproclaim

\demo{Proof} See \cite{LM3}. \qed
\enddemo

\proclaim{Proposition 1.6}
  Let $X$ be a complex projective manifold of dimension $n$ and let
  $\Cal{E}$ be an ample vector bundle of rank $r\geq 2$ on $X$.
  If $g(X,\det\ce)=2$,
  then $n=2$ and $(X,\ce)$ is one of the following:
\roster
\item"(1)"
      $X$ is the Jacobian variety of a smooth curve $B$ of genus $2$
      and $\ce\simeq\ce_r(B,o)\otimes N$ for some $N\in\pic X$
      with $N\equiv 0$, where $\ce_r(B,o)$ is the Jacobian bundle
      for some point $o$ on $B$;
\item"(2)"
      $X\simeq\Bbb{P}_B(\Cal{F})$ for some stable vector bundle $\cf$  
      of rank $2$ on an elliptic curve $B$ with $c_1(\cf)=1$.
      There is an exact sequence 
      $$ 0\to\co_X[2\hf+\rho^*G]\to\ce\to\co_X[\hf+\rho^*T]\to 0, $$
      where $G,T\in\pic B$ and $\rho$ is the projection $X\to B$.
      We have $(\deg G,\deg T)\mathbreak
                              =(-2,1)$ or $(-1,0)$;
\item"($2^{\sharp}$)"
      $X,\cf,B$ and $\rho$ are as in {\rm(2)} 
      and $\Cal{E}\simeq\rho^{*}\Cal{G}\otimes H(\Cal{F})$
      for some stable vector bundle $\cg$ of rank $3$ on $B$
      with $c_{1}(\Cal{G})=-1$;
\item"(3)"
      $X\simeq\Bbb{P}_B(\Cal{F})$ and
      $\Cal{E}\simeq\rho^{*}\Cal{G}\otimes H(\Cal{F})$
      for some semistable vector bundles $\cf$ and $\cg$  
      of rank $2$ on an elliptic curve $B$ with 
      $(c_1(\cf),c_1(\cg))=(1,0)$ or $(0,1)$;
\item"(4)"
      $-K_X$ is ample, $K_X^2=1$ and $\det\ce=-2K_X$.
      We have $\Cal{E}\simeq[-K_X]^{\oplus 2}$, or $c_2(\ce)=3$ and $r=2$;
\item"($5_0$)"
      $X\simeq\Bbb{P}^{1}\times\Bbb{P}^{1}$ and 
      $\ce\simeq\Cal{O}(1,1)\oplus\Cal{O}(1,2)$;
\item"($5_1$)"
      $X$ is the blowing-up of $\Bbb{P}^{2}$ at a point
      and $\Cal{E}\simeq[2L-E]^{\oplus 2}$, where $L$ is the pull-back of
      $\Cal{O}_{\Bbb{P}^{2}}(1)$ and $E$ is the exceptional curve.
\endroster
\endproclaim

\demo{Proof}
See (2.25) Theorem in \cite{Fj2}. \qed
\enddemo

\proclaim{Proposition 1.7}
  Let $X$ be a complex projective manifold of dimension $n$ and let
  $\Cal{E}$ be an ample and spanned vector bundle of rank $r\geq 2$ on $X$.
  If $g(X,\det\ce)=3$,
  then $n=2$ and $(X,\ce)$ is one of the following:
\roster
\item"(1a)"
  $X=\Bbb{P}^{2}$, $\Cal{E}=\Cal{O}_{\Bbb{P}^{2}}(1)^{\oplus 4}$;
\item"(1b)"
  $X=\Bbb{P}^{2}$, and either 
  $\Cal{E}=\Cal{O}_{\Bbb{P}^{2}}(1)^{\oplus 2}\oplus
           \Cal{O}_{\Bbb{P}^{2}}(2)$ 
  or $\Cal{E}=T_{\Bbb{P}^{2}}\oplus\Cal{O}_{\Bbb{P}^{2}}(1)$;
\item"(1c)"
  $X=\Bbb{P}^{2}$, $\operatorname{rank}\Cal{E}=2$ and 
  $\det\Cal{E}=\Cal{O}_{\Bbb{P}^{2}}(4)$;
\item"(2a)"
  $X=\ff_{0}$, and either $\Cal{E}=[\sigma+f]\oplus[\sigma+3f]$
  or $\Cal{E}=[\sigma+2f]^{\oplus 2}$;
\item"(2b)"
  $X=\ff_{1}$, $\Cal{E}=[\sigma+2f]\oplus[\sigma+3f]$;
\item"(2c)"
  $X=\ff_{2}$, $\Cal{E}=[\sigma+3f]^{\oplus 2}$;
\item"(3)"
  $X$ is a Del Pezzo surface with $K_{X}^{2}=2$ and either 
  $\Cal{E}=[-K_{X}]^{\oplus 2}$, or $\Cal{E}=\psi^{*}(\Cal{Q}|_{Y})$,
  where $\psi$ is a birational morphism from $X$ to a surface $Y$ 
  of bidegree $(4,4)$ in the Grassmannian of lines of $\Bbb{P}^{3}$,
  and $\Cal{Q}$ is the universal rank $2$
  quotient bundle;
\item"(4)"
  $X=\Bbb{P}(\Cal{F})$, where $\Cal{F}$ is a rank $2$ vector bundle  
  on an elliptic curve $B$ with $c_{1}(\Cal{F})=1$ and 
  $\Cal{E}=H(\Cal{F})\otimes\rho^{*}\Cal{G}$,
  where $\rho :X\to B$ is the bundle projection and 
  $\Cal{G}$ is any rank $2$ vector bundle on $B$ defined by 
  a non splitting exact sequence
  $0\to\Cal{O}_{B}\to \Cal{G}\to\Cal{O}_{B}(x)\to 0$,
  where $x\in B$.
\endroster
\endproclaim

\demo{Proof} See (1.10) Theorem in \cite{BiLL}. \qed
\enddemo

\comment

\proclaim{Proposition 1.8}
  Let $X$ be a complex projective manifold of dimension $n$ and let
  $\Cal{E}$ be an ample vector bundle of rank $n-2\geq 2$ on $X$
  If $K_X+\det\ce$ is not nef,
  then $(X,\ce)$ is one of the following:
  \roster
  \item $(\pp^n,\co(1)^{\oplus(n-2)})$;
  \item $(\pp^n,\co(1)^{\oplus(n-3)}\oplus\co(2))$;
  \item $(\pp^n,\co(1)^{\oplus(n-3)}\oplus\co(3))$;
  \item $(\pp^n,\co(1)^{\oplus(n-4)}\oplus\co(2)^{\oplus2})$;
  \item $(\qq^n,\co(1)^{\oplus(n-2)})$;
  \item $(\qq^n,\co(1)^{\oplus(n-3)}\oplus\co(2))$;
  \item $(\qq^4,\cs(2))$, where $\cs$ is the spinor bundle 
               on $\qq^4$; 
  \item There is an effective divisor $E$ on $X$ such
                that $(E,\ce_E)\simeq(\pp^{n-1},\co(1)^{\oplus(n-2)})$
                and $\co_E(E)=\co_{\pp^{n-1}}(-1)$;
  \item $X$ is a Fano manifold of index $n-1$ with
                $\pic X\simeq\zz$, and $\ce=H^{\oplus(n-2)}$
                for the ample generator $H$ of $\pic X$;
  \item $X\simeq\pp_B(\cf)$ for some vector bundle $\cf$ 
               of rank $n$ on a smooth curve $B$, and
               $\ce_F=\co_{\pp^{n-1}}(1)^{\oplus(n-2)}$ 
               for every fiber $F$ of the projection $\pi:X\to B$;
  \item $X\simeq\pp_B(\cf)$ for some vector bundle $\cf$ 
                of rank $n$ on a smooth curve $B$, and
                $\ce_F=\co_{\pp^{n-1}}(1)^{\oplus(n-3)}\oplus
                       \co_{\pp^{n-1}}(2)$ 
                for every fiber $F$ of the projection $X\to B$;
  \item There exists a surjective morphism $f:X\to B$ 
                 over a smooth curve $B$ such that
                 $\ce_F=\co_{\qq^{n-1}}(1)^{\oplus(n-2)}$ for 
                 a general fiber $F\simeq\qq^{n-1}$ of $f$;
  \item $X\simeq\pp_S(\cf)$ for some vector bundle $\cf$
                  of rank $n-1$ on a smooth surface $S$, and
                  $\ce_F=\co_{\pp^{n-2}}(1)^{\oplus(n-2)}$ 
                  for every fiber $F$ of the projection $X\to S$.
  \endroster
\endproclaim

\demo{Proof} See Theorem  in \cite{M}. \qed
\enddemo

\endcomment

\subhead
  2. The cases $g(X,\ce)=2$ and $g(X,\Cal{E})=3$
\endsubhead

\proclaim{Theorem 2.1}
  Let $X$ be a compact complex manifold of dimension $n$
  and $\ce$ an ample vector bundle of rank $r$ on $X$ with 
  $1<r<n-1$ and the property $(\ast)$ in {\rm(1.2)}.
  If $g(X,\ce)=2$, then $(X,\ce)$ is one of the following:
  \roster
  \item"(i)" There exists an ample line bundle $A$ on $X$ such that
             $(X,A)$ is a Del Pezzo $4$-fold of degree $1$  
             and $\ce=A^{\oplus2}$ (see also {\rm(2.2.1)});
  \item"(ii)" $X\simeq\pp_B(\cf)$ and $\ce=\hf\otimes\pi^*\cg$,
              where $\cf$ and $\cg$ are vector bundles on
              an elliptic curve $B$ such that $\rk\cf=4$, $\rk\cg=2$,
              $c_1(\cf)+2c_1(\cg)=1$, and $\pi:X\to B$ is the
              bundle projection;
  \item"(iii)" $X\simeq\pp_B(\cf)$ and $\ce=\hf\otimes\pi^*\cg$,
               where $\cf$ and $\cg$ are vector bundles on
               an elliptic curve $B$ such that $\rk\cf=5$, $\rk\cg=3$,
               $3c_1(\cf)+5c_1(\cg)=1$, and $\pi:X\to B$ is the
               bundle projection.
  \endroster
\endproclaim

\demo{Proof}
  Suppose that $g(X,\ce)=2$.
  Since $\ce$ satisfies $(\ast)$, there exists a nonzero section
  $s\in H^0(X,\ce)$ whose zero locus $Z:=(s)_0$ is a smooth submanifold
  of $X$ of dimension $n-r$ and $2=g(X,\ce)=g(Z,\det\ce_Z)$.
  From (1.6) we see that $n-r=2$ and $(Z,\ce_Z)$ is one of the cases
  in (1.6).
  We make a case by case analysis in the following.
 
  (2.1.1) If $(Z,\ce_Z)$ is in case (1.6;1), then $K_Z=\co_Z$.
  We have $K_X+\det\ce=\co_X$ since $[K_X+\det\ce]_Z=K_Z$ and
  $\pic(\iota):\pic(X)\to\pic(Z)$ is injective by (1.3).
  We get also that $H^1(\iota):H^1(X,\zz)\to H^1(Z,\zz)$ is
  an isomorphism by (1.3), but this is impossible since
  $X$ is a Fano manifold and $Z$ is an abelian surface.  

  (2.1.2) If $(Z,\ce_Z)$ is in case (1.6;$5_0$), we have $r=2$ and $n=4$.
  By (1.4), $(X,\ce)$ is one of the cases (Q1),(Q2) and (Q3).
  We easily see that $g(X,\ce)\not=2$ in cases (Q1) and (Q2).  
  In case (Q3), we can write $\cf=\oplus_{i=1}^4\co_{\pp^1}(a_i)$. 
  Since $\ce$ is ample, $\hf+\pi^*\co_{\pp^1}(b_j)$ is ample 
  and so is $\cf\otimes\co_{\pp^1}(b_j)$.
  Hence we get $a_i+b_j>0$ for every $i$ and $j$.
  Then it follows that 
  $$
    2=2g(X,\ce)-2=(K_X+2c_1(\ce))c_1(\ce)c_2(\ce)
     =2(-2+\sum_{i=1}^4 a_i+2(b_1+b_2))\ge 4,
  $$ 
  a contradiction.

  (2.1.3) If $(Z,\ce_Z)$ is in case (1.6;$5_1$), we have $r=2$ and $n=4$.
  Since $Z=\ff_1$, we see that $(X,\ce)$ is in case (1.5;3).
  If $(X,\ce)$ is the type (1.5;1) with $B=\pp^1$, then we come to
  a contradiction by the argument of (2.1.2).   
  Hence we have 
  $X\simeq \Bbb{P}_{\Bbb{P}^{2}}(\Cal{F})$ for some ample vector bundle
  $\Cal{F}$ on $\Bbb{P}^{2}$ with $c_{1}(\Cal{F})=2k+3~(k>0)$, and 
  $\Cal{E}_{F}=\Cal{O}_{\Bbb{P}^2}(1)^{\oplus2}$
  for every fiber $F$ of the bundle map $\pi: X\to \pp^2$.
  We set $H:=\hf$;
  we can write $\ce=H\otimes\pi^*\cg$ 
  for some vector bundle $\cg$ of rank $2$ on $\pp^2$.
  Since $\ce_Z=[2L-E]^{\oplus2}$, we can write $H_Z=aL-E~(2\le a\in\zz)$.
  Then we get $\cg=\co_{\pp^2}(2-a)^{\oplus2}$,
  hence $\ce=[H+\pi^*\co_{\pp^2}(2-a)]^{\oplus2}$
  by $(\pi|_Z)^*\cg=\ce_Z\otimes[-H_Z]=[(2-a)L]^{\oplus2}$.
  Since $\ce$ is ample, $H+\pi^*\co_{\pp^2}(a)$ is ample and 
  so is $\cf\otimes\co_{\pp^2}(a)$.
  Then we get $c_1(\cf\otimes\co_{\pp^2}(2-a))\ge 3$,
  hence $2k-3a+6\ge 0$.
  We note that 
  $$ 3=(2L-E)^2=c_2(\ce_Z)=c_2(\ce)^2
      =s_2(\cf)+4c_1(\cf)\cdot(2-a)+6(2-a)^2.
  $$
  On the other hand, we have
  $$ a^2-1=(aL-E)^2=H_Z^2=H^2\cdot c_2(\ce)
          =s_2(\cf)+2c_1(\cf)\cdot(2-a)+(2-a)^2.
  $$
  From these two equalities we get $(2-a)(2k-3a+7)=0$.
  Since $2k-3a+6\ge 0$, we have $a=2$ and then $c_2(\cf)=3$ 
  and $\ce=H^{\oplus2}$.
  It follows that 
  $$ 2=2g(X,\ce)-2=(K_X+2c_1(\ce))c_1(\ce)c_2(\ce)
      =2s_2(\cf)+4k\ge 10,
  $$
  a contradiction.

  (2.1.4) If $(Z,\ce_Z)$ is in case (1.6;4), then $r=2$ and $n=4$.
  We have $2K_X+3\det\ce=\co_X$ since $[2K_X+3\det\ce]_Z=\co_Z$
  and the restriction map $\pic(X)\to\pic(Z)$ is injective.
  By setting $A:=K_X+2\det\ce$, we get $\det\ce=2A$ and $K_X+3A=\co_X$,
  hence $(X,A)$ is a Del Pezzo 4-fold.
  Then we set $\ce':=\ce\oplus A$; we get $K_X+\det\ce'=\co_X$ and
  $\ce'\simeq A^{\oplus3}$ by using \cite{PSW,(7.4) Proposition}.
  It follows that $\ce\simeq A^{\oplus2}$ and
  $$ 2=2g(X,\ce)-2=(K_X+2c_1(\ce))c_1(\ce)c_2(\ce)=2A^4,
  $$
  hence $A^4=1$.
  Thus we obtain that $(X,\ce)$ is the case (i) of our theorem.
 
  (2.1.5) If $(Z,\ce_Z)$ is in case (1.6;2), then $r=2$ and $n=4$.
  Since $Z$ is a geometrically ruled surface over an elliptic curve $B$,
  by (1.5), $X$ is a $\pp^3$-bundle over $B$ and 
  $\ce|_{\widetilde{F}}=\co_{\pp^3}(1)^{\oplus2}$ for every fiber 
  $\widetilde{F}$ of the ruling $\pi:X\to B$.
  On the other hand, we have $\ce_Z|_F=\co_{\pp^1}(1)\oplus\co_{\pp^1}(2)$ 
  for every fiber $F$ of the ruling $\rho:Z\to B$. 
  This is a contradiction since $\pi|_Z=\rho$.
  If $(Z,\ce_Z)$ is in case (1.6;$2^{\sharp}$) or (1.6;3),
  by using (1.5), we obtain that $(X,\ce)$ is the case (ii) or (iii) 
  of our theorem respectively.
  This completes the proof. \qed
\enddemo

\example{Remark 2.2}
 We make some comments on (2.1).

 (2.2.1) In case (2.1; i), Del Pezzo $4$-folds of degree $1$
  have been classified in \cite{Fj1, Part III}.
  In particular, they are weighted hypersurfaces of degree $6$ in 
  the weighted projective space $\pp(3,2,1,1,1,1)$.

 (2.2.2)
  We give an example of $(X,\ce)$ in case (2.1; ii) in the following.
  Let $L_1$ and $L_2$ be line bundles on an elliptic curve $B$
  such that $\deg L_1=\deg L_2$ and $L_1\not= L_2$ in $\pic B$.
  Let $\cf$ be an indecomposable vector bundle of rank 4 on $B$
  with $c_1(\cf)=1-2\deg L_1-2\deg L_2$.
  We set $X:=\pp_B(\cf)$, $\cg:=L_1\oplus L_2$, and
  $\ce:=\hf\otimes\pi^*\cg=\oplus_{i=1}^2[\hf+\pi^*L_i]$,
  where $\pi:X\to B$ is the bundle projection.
  Since $c_1(\cf\otimes L_i)=1$, $\cf\otimes L_i$ is ample and
  $h^0(B,\cf\otimes L_i)=1$. 
  Then there exists an exact sequence
  $$ 0\to\co_B\to\cf\otimes L_i\to Q_i\to 0,
  $$
  where $Q_i$ is a locally free sheaf of rank 3 on $B$.
  Since $Q_i$ is ample and $c_1(Q_i)=1$, we see that $Q_i$ is 
  indecomposable.
  We set $D_i:=\pp_B(Q_i)$ and $Z:=D_1\cap D_2$.
  Since $c_1(Q_2\otimes[L_1-L_2])=1$, there exists an exact sequence
  $$ 0\to\co_B\to Q_2\otimes[L_1-L_2]\to Q\to 0,
  $$ 
  where $Q$ is a locally free sheaf of rank 2 on $B$.
  Then we have $Z=\pp_B(Q)$ in $|H(Q_2)+(\pi|_{D_2})^*(L_1-L_2)|$.
  Thus we see that $(X,\ce)$ satisfies the condition $(\ast)$
  and $(X,\ce)$ is an example of (2.1; ii).

 (2.2.3) The authors have no example for case (2.1; iii).
  We note that without the condition $(\ast)$ 
  we have examples for the case.
  Indeed, we can take semistable vector bundles $\cf$ and $\cg$
  on an elliptic curve $B$ with the property that
  $\rk\cf=5$, $\rk\cg=3$, and $3c_1(\cf)+5c_1(\cg)=1$.
  Let $\pi:\pp(\cf)\to B$ and $\pi':\pp(\cg)\to B$ be the 
  bundle projections.
  Then $5\hf-\pi^*\det\cf$ is nef on $\pf$ and
  $3\hg-(\pi')^*\det\cg$ is nef on $\pg$ by \cite{Mi, Theorem 3.1}.
  We set $\ce:=\hf\otimes\pi^*\cg$ and let $p:\pe\to B$ be the 
  composition of the projection $\pe\to\pf$ and $\pi$.
  Then $15\he-F$ is nef on $\pe$ for a fiber $F$ of $p$,
  hence $\ce$ is ample.
  But it is uncertain that such $\ce$ satisfies $(\ast)$.

 (2.2.4) We see that every vector bundle $\ce$ appeared in (2.1)
  is not spanned.
  Indeed, it is clear for case (2.1; i).
  For cases (2.1; ii) and (2.1; iii), we use the following
\proclaim{Lemma 2.2.5}
  Let $\cf$ be a vector bundle of rank $r$ on an elliptic curve.
  Then there exists a line sub-bundle $L$ of $\cf$ such that
  $\deg L\ge[c_1(\cf)/r]$, where $[c_1(\cf)/r]$ is the largest integer
  that is not greater than $c_1(\cf)/r$.
\endproclaim
  This is a consequence of the Mukai-Sakai theorem \cite{MuS},
  hence proof is omitted.

  Suppose that $\ce$ is spanned in case (2.1; ii).  
  Applying the lemma to $\cf^{\vee}$ and $\cg^{\vee}$,
  we get quotient line bundles $L_1$ and $L_2$ of $\cf$ and $\cg$
  respectively,
  with the property that $\deg L_1\le-[-c_1(\cf)/4]$ and
  $\deg L_2\le-[-c_1(\cg)/2]$.
  The surjection $\cf\to L_1$ gives a section $C:=\pp(L_1)$
  of the projection $\pi:\pp_B(\cf)\to B$.
  Since $\hf|_C=(\pi|_C)^*L_1$, we see that $(\pi|_C)^*(L_1\otimes L_2)$ 
  is a quotient line bundle of $\ce_C$,
  hence $L_1\otimes L_2$ is ample and spanned.
  From $c_1(\cf)+2c_1(\cg)=1$ we get 
  $\deg L_1+\deg L_2\le-[(2c_1(\cg)-1)/4]-[-c_1(\cg)/2]=1$;
  this leads to a contradiction since $B$ is an elliptic curve.
  Similarly we can show that $\ce$ is not spanned in case (2.1; iii).
\endexample

\proclaim{Theorem 2.3}
  Let $X$ be a compact complex manifold of dimension $n$
  and $\ce$ an ample and spanned vector bundle of rank $r$ on $X$ 
  with $1<r<n-1$.
  If $g(X,\ce)=3$, then $(X,\ce)$ is one of the following:
  \roster
  \item"(i)" $(\pp^6,\co_{\pp^6}(1)^{\oplus4})$;
  \item"(ii)" $(\pp^1\times\pp^3,\co_{\pp^1\times\pp^3}(1,1)^{\oplus2})$;
  \item"(iii)" There exists a double covering $f:X\to\pp^4$
               with branch locus $B\in|\co_{\pp^4}(4)|$ and
               $\ce=f^*\co_{\pp^4}(1)^{\oplus2}$.
  \endroster
\endproclaim

\demo{Proof}
  Suppose that $g(X,\ce)=3$.
  We argue as in the proof of (2.1).
  Since $\ce$ is spanned, there exists a nonzero section
  $s\in H^0(X,\ce)$ whose zero locus $Z:=(s)_0$ is a smooth submanifold
  of $X$ of dimension $n-r$ and $3=g(X,\ce)=g(Z,\det\ce_Z)$.
  From (1.7) we see that $n-r=2$ and $(Z,\ce_Z)$ is one of the cases
  in (1.7).
 
  (2.3.1) If $(Z,\ce_Z)$ is in case (1a),(1b), or (1c) of (1.7), 
  then $Z=\pp^2$ and $(X,\ce)$ is the case (P1) of (1.4) since $n-r=2$.
  We obtain that $(X,\ce)$ is the case (i) of our theorem
  by $g(X,\ce)=3$.

  (2.3.2) If $(Z,\ce_Z)$ is in case (3) of (1.7), then $r=2$ and $n=4$.
  By setting $A:=K_X+2\det\ce$, we infer that $(X,A)$ is a Del Pezzo
  manifold and $\ce=A^{\oplus2}$ from the same argument as that in (2.1.4).
  Then we find that $A^4=2$ since $g(X,\ce)=3$.
  Hence we obtain that $(X,\ce)$ is the case (iii) of our theorem
  by \cite{Fj1, Part I}.
  
  (2.3.3) If $(Z,\ce_Z)$ is in case (2a),(2b),(2c), or (4) of (1.7), 
  then $r=2$ and $n=4$.
  Since $Z$ is a geometically ruled surface, by (1.5),
  $(X,\ce)$ is one of the following:
  \roster
  \item"(R1)" $(\Bbb{P}^{4}, \Cal{O}_{\Bbb{P}^{4}}(1)\oplus
                               \Cal{O}_{\Bbb{P}^{4}}(2))$;
  \item"(R2)" $(\Bbb{Q}^{4}, 
                \Cal{O}_{\Bbb{Q}^{4}}(1)^{\oplus2})$;
  \item"(R3)" $X$ is a $\Bbb{P}^3$-bundle over a smooth curve $B$ and 
              $\Cal{E}_{\widetilde{F}}=\Cal{O}_{\Bbb{P}^{3}}(1)^{\oplus2}$
              for every fiber $\widetilde{F}$ of the bundle map 
              $\pi: X\to B$;
  \item"(R4)" $Z=\ff_1$, $X\simeq \Bbb{P}_{\Bbb{P}^{2}}(\Cal{F})$ 
              for some ample vector bundle
              $\Cal{F}$ on $\Bbb{P}^{2}$ with 
              $c_{1}(\Cal{F})=2k+3$ ($k>0$), and 
              $\Cal{E}_{\widetilde{F}}=\Cal{O}_{\Bbb{P}^{2}}(1)^{\oplus2}$
              for every fiber $\widetilde{F}$ of the bundle map 
              $\pi: X\to \pp^2$.
  \endroster 
  Cases (R1) and (R2) are ruled out by $g(X,\ce)=3$.
  Case (R4) comes from (2b) of (1.7), hence $\pi|_Z$ is the blowing-up 
  $\ff_1\to\pp^2$ and $\Cal{E}_Z=[\sigma+2f]\oplus[\sigma+3f]$.
  We can write $\ce=\hf\otimes\pi^*\cg$ for some vector bundle $\cg$
  of rank 2 on $\pp^2$ and $\hf_Z=a\sigma+bf$ for some $a,b\in\zz$.
  Then 
  $$ 2\sigma+5f=\det\ce_Z=2\hf_Z+(\pi|_Z)^*\det\cg
               =(2a+c_1(\cg))\sigma+(2b+c_1(\cg))f,
  $$
  hence $2a-2b=-3$, a contradiction.
  In case (R3), we have $X\simeq\pp_B(\cf)$ and 
  $\ce=\hf\otimes\pi^*\cg$ for some vector bundles $\cf$ and $\cg$
  on $B$ such that $\rk\cf=4$ and $\rk\cg=2$.
  Then
  $$ 4=2g(X,\ce)-2=(K_X+2c_1(\ce))c_1(\ce)c_2(\ce)
                  =2(2g(B)-2+c_1(\cf)+2c_1(\cg)), 
  $$
  where $g(B)$ is the genus of $B$.
  Since $\ce$ is ample, we find that $c_1(\cf)+2c_1(\cg)>0$
  from $(\det\ce)^4>0$.
  It follows that $g(B)\le 1$.
  In case $g(B)=0$, we have $B\simeq\pp^1$ and $c_1(\cf)+2c_1(\cg)=4$.
  Then we can write $\cf=\sum_{i=1}^4\co(a_i)$ and
  $\cg=\sum_{j=1}^2\co(b_j)$.
  By the same argument as that in (2.1.2), we infer that 
  $a_i+b_j=1$ for every $i$ and $j$.
  It follows that $a_1=\dots=a_4$ and $b_1=b_2$,
  hence $\pp_B(\cf)\simeq\pp^1\times\pp^3$ and
  $\ce=\co_{\pp^1\times\pp^3}(1,1)^{\oplus2}$,
  which is the case (ii) of our theorem.
  In case $g(B)=1$, we have $c_1(\cf)+2c_1(\cg)=2$.
  Then we get a contradiction by the same argument as that in (2.2.4).
  We have thus completed the proof. \qed
\enddemo

\subhead
  3. The cases $g(X,\Cal{E})=q(X)+1$ and $g(X,\Cal{E})=q(X)+2$
\endsubhead

\proclaim{Theorem 3.1}
Let $X$ be a compact complex manifold of dimension $n$
and let $\Cal{E}$ be an ample and spanned vector bundle of rank $r$
with $1<r<n-1$.
Then $g(X,\Cal{E})=q(X)+1$ if and only if
$(X,\Cal{E})$ is one of the following:
  \roster
  \item $(\pp^5,\co_{\pp^5}(1)^{\oplus2})$;
  \item $(\pp^5,\co_{\pp^5}(1)^{\oplus3})$;
  \item $(\qq^4,\co_{\qq^4}(1)^{\oplus2})$.
  \endroster
\endproclaim

\demo{Proof}
First we note that if $(X,\ce)$ is one of the cases (1),(2) and (3)
of our theorem,
then we easily see that $g(X,\Cal{E})=1=q(X)+1$.
Suppose that $g(X,\Cal{E})=q(X)+1$ on the contrary.
Let $Z$ be a smooth submanifold of $X$ with $\dim Z=n-r$
defined as the zero locus of some $s\in H^{0}(X,\Cal{E})$.
Then $g(X,\Cal{E})=g(Z, \operatorname{det}\Cal{E}_{Z})$.
We put $A:=\operatorname{det}\Cal{E}_{Z}$;
then $A$ is ample and spanned.
If $n-r\ge 3$,
we take general members $D_1,\dots,D_{n-r-2}\in|A|$
with the property that $S:=D_1\cap\dots\cap D_{n-r-2}$
is a smooth surface.
If $n-r=2$, we set $S=Z$.
Hence there exists a polarized surface $(S,A_{S})$
such that $g(Z,A)=g(S,A_{S})$.
We get $q(X)=q(Z)=q(S)$ by using (1.3).
Thus we get $g(S,A_{S})=q(S)+1$.

(I) The case in which $\kappa(S)=-\infty$.
    ($\kappa(S)$ is the Kodaira dimension of $S$.)

First we remark that $(S,A_{S})$ is not a scroll over a smooth curve
by construction.

\proclaim{Claim 3.2}
 $g(S,A_{S})\geq 2q(S)$.
\endproclaim

\demo{Proof}
If $q(S)=0$, then this is obvious.
So we may assume that $q(S)\geq 1$.
Then by the above remark, $K_{S}+A_{S}$ is nef.
So we get
$$\split
0\leq (K_{S}+A_{S})^{2}&=(K_{S})^{2}+2(K_{S}+A_{S})A_{S}-(A_{S})^{2} \\
                       &\leq 8(1-q(S))+4(g(S,A_{S})-1)-(A_{S})^{2} \\
                       &=4(g(S,A_{S})-2q(S)+1)-(A_{S})^{2}.
\endsplit
$$
Hence $g(S,A_{S})\geq 2q(S)$. \qed
\enddemo

By (3.2) and the above argument, we get $q(S)\leq 1$ and 
$g(X,\Cal{E})\leq 2$.
So we obtain that $(X,\ce)$ is the case (1),(2), or (3) 
of our theorem by using (2.1) and \cite{I}.

(II) The case in which $\kappa(S)\geq 0$.

By Riemann-Roch Theorem and Vanishing Theorem, we get 
 $$h^{0}(K_{S}+A_{S})-h^{0}(K_{S})=g(S,A_{S})-q(S)=1.$$
If $h^{0}(K_{S})>0$, then 
 $$h^{0}(K_{S}+A_{S})\geq h^{0}(K_{S})+h^{0}(A_{S})-1.$$
But this is impossible since $h^{0}(A_{S})\geq 3$.
Hence $h^{0}(K_{S})=0$.
By the classification theory of surfaces, we get $q(S)\leq 1$.
Hence $g(X,\Cal{E})\leq 2$.
From (2.1) and \cite{I} we infer this case cannot occur.
This completes the proof. \qed
\enddemo

\example{Remark 3.3}
Let $L$ be an ample and spanned line bundle on a compact complex 
manifold $X$ of dimension $n\ge 2$.
When $n\ge 3$, we have $g(X,L)=q(X)+1$ if and only if 
$(X,L)$ is a Del Pezzo manifold (see \cite{Fk3}).
When $n=2$, we have $g(X,L)=q(X)+1$ if and only if 
$(X,L)$ is a Del Pezzo surface (i.e. $L=-K_X$) or
$X\simeq\pp_B(\cf)$ and $L\equiv 2\hf$ for some ample vector bundle
$\cf$ of rank $2$ on an elliptic curve $B$ with $c_1(\cf)=1$.
We can prove this by the argument in (3.1) and a classification result
\cite{LP, (3.1)}.
\endexample

\proclaim{Proposition 3.4}
Let $X$ be a compact complex manifold of dimension $n$
and let $\Cal{E}$ be an ample and spanned vector bundle of rank $r$
with $1<r<n-1$.
Then we have $g(X,\Cal{E})\not=q(X)+2$.
\endproclaim

\demo{Proof}
The following argument is similar to the proof of (3.1).
Suppose that $g(X,\Cal{E})=q(X)+2$.
Let $Z$ be a smooth submanifold of $X$ with $\dim Z=n-r$
defined as the zero locus of some $s\in H^{0}(X,\Cal{E})$.
Then $g(X,\Cal{E})=g(Z, \operatorname{det}\Cal{E}_{Z})$ and
$\operatorname{det}\Cal{E}_{Z}$ is ample and spanned.
As in the proof of (3.1), we get a smooth surface $S$
such that $g(Z,\det\ce_Z)=g(S,\det\ce_{S})$.
We have $q(X)=q(Z)=q(S)$, thus we get $g(S,\det\ce_{S})=q(S)+2$.
Then we find that $q(S)\le 1$ by \cite{R, Theorem 3.4}. 
It follows that $g(X,\Cal{E})\le 3$ and we infer that 
$(X,\ce)$ does not exist from (2.1) and (2.3). 
This completes the proof. \qed
\enddemo

\example{Remark 3.5}
We see that the pairs $(X,\ce)$ in (2.3) satisfy $g(X,\Cal{E})=q(X)+3$.  
In Appendix we give a classification of polarized surfaces $(X,L)$
such that $g(X,L)=q(X)+2$ and $L$ is spanned.
\endexample

\subhead
  4. Another Lower bound for $g(X,\ce)$
\endsubhead

\proclaim{Proposition 4.1}
Let $L$ be an ample and spanned line bundle on a compact
complex manifold $X$ with $\dim X=n\geq 2$.
Then $g(X,L)\geq 2q(X)-1$ unless $(X,L)$ is a scroll 
over a smooth curve $B$ of genus $g(B)\geq 2$.
\endproclaim

\demo{Proof}
Since $L$ is ample and spanned, if $n\geq 3$, we can take
general members $D_{1},\dots,D_{n-2}\in |L|$ such that
$S:=D_{1}\cap\dots \cap D_{n-2}$ is a smooth surface.
If $n=2$, we set $S=X$.
Then we get $g(X,L)=g(S,L_{S})$ and $q(X)=q(S)$.

If $\kappa(S)\geq 0$, then $g(X,L)=g(S,L_{S})\geq 2q(S)-1=2q(X)-1$ by 
\cite{Fk1, Corollary 3.2}.

If $\kappa(S)=-\infty$ and $(S,L_{S})$ is not a scroll over a smooth curve,
then $g(X,L)=g(S,L_{S})\geq 2q(S)=2q(X)$ by (3.2).

If $\kappa(S)=-\infty$ and $(S,L_{S})$ is a scroll over a smooth curve,
then $g(X,L)=g(S,L_{S})=q(S)=q(X)$.
Hence we get $g(X,L)\geq 2q(X)-1$ if $q(S)\leq 1$.
So we may assume that $q(S)\geq 2$.
Then we obtain that $(X,L)$ is a scroll over a smooth curve $B$
of genus $g(B)\geq 2$ by using Theorem 3 in \cite{B\v a}. \qed
\enddemo

\proclaim{Theorem 4.2}
Let $X$ be a compact complex manifold with $\dim X=n$
and let $\Cal{E}$ be an ample and spanned vector bundle of rank $r$
with $1<r<n-1$.
Then $g(X,\Cal{E})\geq 2q(X)-1$.
\endproclaim

\demo{Proof}
Let $Z$ be the zero locus of some $s\in H^{0}(X,\Cal{E})$
such that $Z$ is a smooth submanifold of $X$ with $\dim Z=n-r$.
Then $g(X,\Cal{E})=g(Z, \operatorname{det}\Cal{E}_{Z})$
and $q(X)=q(Z)$.
We put $A:=\operatorname{det}\Cal{E}_{Z}$;
then $A$ is ample and spanned.
Since $(Z,A)$ is not a scroll, by (4.1), we obtain that
$g(X,\Cal{E})=g(Z,A)\geq 2q(Z)-1=2q(X)-1$. \qed
\enddemo

\subhead
  5. The case of a fiber space over a curve
\endsubhead

\definition{Definition 5.1}
  Here we say that a quartet $(f,X,Y,\Cal{E})$ is a {\it generalized 
  polarized fiber space} if 
  \roster
  \item $X$ and $Y$ are compact complex manifolds with 
        $1\le\dim Y<\dim X$,
  \item $f:X\to Y$ is a surjective morphism with connected fibers, and
  \item $\ce$ is an ample vector bundle on $X$.
  \endroster
\enddefinition

\proclaim{Theorem 5.2}
Let $(f,X,C,\Cal{E})$ be a generalized polarized fiber space with
$n:=\dim X\geq 2$, $\dim C=1$, and $r:=\rk\ce\leq n-1$.
Then $g(X,\ce)\geq g(C)$.
\endproclaim

\demo{Proof}
First we remark that the following equality holds:
$$
\split
g(X,\ce)&=g(C)+\frac{1}{2}(K_{X/C}+(n-r)c_{1}(\ce))
c_{1}(\ce)^{n-r-1}c_{r}(\ce) \\
&\ \ \ +(g(C)-1)(c_{1}(\ce)^{n-r-1}c_{r}(\ce)F-1),
\endsplit$$
where $K_{X/C}:=K_{X}-f^{*}(K_{C})$ and $F$ is a general fiber of $f$.
\flushpar
If $g(C)=0$, then Theorem 5.2 is true by \cite{I}.
So we may assume that $g(C)\geq 1$.
\flushpar
(I) The case in which $K_{X/C}+(n-r)c_{1}(\ce)$ is $f$-nef.
\flushpar
Then there exists a surjective map
$$f^{*}\circ f_{*}(\Cal{O}(m(K_{X/C}+(n-r)c_{1}(\ce))))\to
\Cal{O}(m(K_{X/C}+(n-r)c_{1}(\ce)))$$
for any large $m$ by base point free theorem.

By Theorem A in Appendix in \cite{Fk2},
$f_{*}(\Cal{O}(m(K_{X/C}+(n-r)c_{1}(\ce))))$ is semipositive.
Hence $K_{X/C}+(n-r)c_{1}(\ce)$ is nef.
So we get
$$(K_{X/C}+(n-r)c_{1}(\ce))c_{1}(\ce)^{n-r-1}c_{r}(\ce)\geq 0.$$
Hence we obtain $g(X,\ce)\geq g(C)$ because 
$c_{1}(\ce)^{n-r-1}c_{r}(\ce)F\geq 1$.
\flushpar
(II) The case in which $K_{X/C}+(n-r)c_{1}(\ce)$ is not $f$-nef.
\flushpar
Then $K_{X}+(n-r)c_{1}(\ce)$ is not nef.
So by Mori Theory, there exists an extremal rational curve $l$ 
such that $(K_{X}+(n-r)c_{1}(\ce))l<0$.
Hence 
$$n+1\geq -K_{X}l>(n-r)c_{1}(\ce)l\geq (n-r)r\geq n-1.$$
\flushpar 
If $(n-r)r=n$, then $(n,r)=(4,2)$. \flushpar
If $(n-r)r=n-1$, then $r=1$ or $r=n-1$.
\flushpar
(II-1) The case where $(n,r)=(4,2)$.
\flushpar
Then $-K_{X}l=5=n+1$.
So we have $\operatorname{Pic}X\cong \zz$ by \cite{W}.
But this is impossible because $X$ has a nontrivial fibration.
\flushpar
(II-2) The case in which $r=1$.
\flushpar
Then this is true by Theorem 1.2.1 in \cite{Fk2}.
\flushpar
(II-3) The case in which $r=n-1$.
\flushpar
If $n=2$, then $r=1$ and so we may assume that $n\geq 3$.
Since $X$ has a nontrivial fibration, $(X,\ce)$ is the following
type by \cite{YZ}:
there exists a surjective morphism $\pi : X\to B$ such that 
any fiber of $\pi$ is $\pp^{n-1}$ 
and $\ce|_{F_{\pi}}\cong \Cal{O}(1)^{\oplus n-1}$,
where $B$ is a smooth curve and $F_{\pi}$ is a fiber of $\pi$.
\flushpar
Since any fiber of $\pi$ is $\pp^{n-1}$, there exists a morphism
$\delta: B\to C$ such that $f=\delta\circ\pi$.
Because $f$ has connected fibers, $\delta$ is an isomorphism.
In particular, $g(B)=g(C)$.
On the other hand, by \cite{Ma}, $g(X,\ce)=g(B)$.
Hence $g(X,\ce)=g(C)$.
This completes the proof of Theorem 5.2. \qed
\enddemo  

\proclaim{Theorem 5.3}
Let $(f,X,C,\ce)$ be a generalized polarized fiber space
with $n\geq r+1$ and $\dim C=1$,
where $n:=\dim X\geq 3$ and $r:=\rk\ce\geq 2$.
If $g(X,\ce)=g(C)$, then
any fiber $F$ of $f$ is isomorphic to $\pp^{n-1}$ and
$\ce|_{F}\cong \oplus_{i=1}^{n-1}\Cal{O}_{\pp^{n-1}}(1)$.
\endproclaim

\demo{Proof}
(I) The case in which $g(C)\leq 1$.
\flushpar
Then $g(X,\ce)=g(C)\leq 1$, and by the classification results of 
\cite{I} and \cite{Ma}, we get the following:
$X$ is a $\pp^{n-1}$-bundle over $\pp^{1}$ or a smooth elliptic curve
and $\Cal{E}|_{F_{\pi}}\cong \Cal{O}_{\pp^{n-1}}(1)^{\oplus n-1}$,
where $F_{\pi}$ is a fiber of its bundle map $\pi: X\to B$ and 
$B$ is $\pp^{1}$ or 
a smooth elliptic curve.
Since any fiber of $\pi$ is $\pp^{n-1}$, there exists a morphism
$\delta: B\to C$ such that $f=\delta\circ\pi$.
Because $f$ has connected fibers, $\delta$ is an isomorphism.
Therefore we get the assertion.
\flushpar
(II) The case in which $g(C)\geq 2$.
\flushpar
(II-1) $n-r\geq 2$ case.
\flushpar
If $K_{X/C}+(n-r-1)c_{1}(\ce)$ is $f$-nef, then by the same argument 
as in the proof of Theorem 5.2 we get
$$(K_{X/C}+(n-r-1)c_{1}(\ce))c_{1}(\ce)^{n-r-1}c_{r}(\ce)\geq 0$$
and 
$$(K_{X/C}+(n-r)c_{1}(\ce))c_{1}(\ce)^{n-r-1}c_{r}(\ce)\geq 1.$$
Hence we obtain that $g(X,\ce)>g(C)$. 
So we may assume that $K_{X/C}+(n-r-1)c_{1}(\ce)$ is not $f$-nef.
Then by Mori Theory, there exists an extremal rational curve $l$ 
such that $(K_{X}+(n-r-1)c_{1}(\ce))l<0$.
Hence we get 
$$n+1\geq -K_{X}l>(n-r-1)c_{1}(\ce)l\geq (n-r-1)r\geq n-2.$$
\flushpar 
If $(n-r-1)r=n$, then $-K_{X}l=n+1$ and $\operatorname{Pic}X\cong\zz$
by \cite{W}.
But this is impossible.
\flushpar
If $(n-r-1)r=n-1$, then $n=5$ and $r=2$.
\flushpar
Here we prove the following Lemma.

\proclaim{Lemma 5.4}
Let $(f,X,C,\ce)$ be a generalized polarized fiber space
with $n\geq r+1$, $\dim C=1$ and $g(C)\geq 1$,
where $n:=\dim X\geq 3$ and $r:=\rk\ce\geq 2$.
If $\kappa(K_{F}+xc_{1}(\ce_{F}))\geq 0$ for a rational number $x$ with
$x<n-r$, and a general fiber $F$ of $f$, then $g(X,\ce)\geq g(C)+1$.
\endproclaim

\demo{Proof}
By assumption, there exists a natural number $N$
such that 
$$f_{*}(\Cal{O}(N(K_{X/C}+xc_{1}(\ce))))\not=0.$$
By Remark 1.3.2 in \cite{Fk2}, 
$N(K_{X/C}+xc_{1}(\Cal{E}))$ is pseudo effective.
Therefore 
$$(K_{X/C}+xc_{1}(\ce))c_{1}(\ce)^{n-r-1}c_{r}(\ce)\geq 0$$
and we get 
$$(K_{X/C}+(n-r)c_{1}(\ce))c_{1}(\ce)^{n-r-1}c_{r}(\ce)\geq 1.$$
Since $g(C)\geq 1$, we get that $g(X,\ce)\geq g(C)+1$. \qed
\enddemo

We continue the proof of Theorem 5.3.
If $K_{F}+xc_{1}(\ce)$ is nef for a rational number $x$
with $x<3$, then we can prove that $g(X,\ce)>g(C)$.
\flushpar
Assume that 
$K_{F}+xc_{1}(\ce_{F})$ is not nef for a rational number $x$ with $x<3$.
Then there exists an extremal rational curve $l$ on $F$ such that
$n\geq -K_{F}l>xc_{1}(\ce_{F})l\geq rx$.
Since $n=5$ and $r=2$, we have $x<5/2$.
Therefore there exists a rational number $y<3$ such that
$K_{F}+yc_{1}(\ce_{F})$ is nef, and we get $g(X,\ce)>g(C)$.
\flushpar
If $(n-r-1)r=n-2$, then $r=n-2$ by assumption.
Assume that $K_{F}+xc_{1}(\ce_{F})$
is not nef for a rational number $x$ with $x<2$.
Then we get $n>rx$ by the same argument as above.
Since $r=n-2$, we get $x<n/(n-2)=1+2/(n-2)$.
By assumption, we get $n\geq 4$.
So we have $x<2$.
Therefore there exists a rational number $y<2$
such that $K_{F}+yc_{1}(\ce_{F})$ is nef.
Hence we have $g(X,\ce)>g(C)$.
\flushpar
(II-2) $n-r=1$ case.
\flushpar
First we assume that $K_{F}+c_{1}(\ce_{F})$ is nef for a general fiber 
$F$ of $f$.
If $K_{F}+c_{1}(\ce_{F})$ is ample, then there exists a rational number
$t>0$ such that $\kappa(K_{F}+(1-t)c_{1}(\ce_{F}))\geq 0$ by Kodaira's
Lemma.
So we get that $g(X,\ce)>g(C)$ by the same argument as above.
Assume that $K_{F}+c_{1}(\ce_{F})$ is not ample.
Since $\dim F=\rk\ce_{F}$, by \cite{Fj3},
we get that $(F,\ce_{F})$ is one of the following:
\roster
\item"(A)" $(\Bbb{P}^{n-1}, \Cal{O}_{\Bbb{P}^{n-1}}(2)\oplus 
                           \Cal{O}_{\Bbb{P}^{n-1}}(1)^{\oplus n-2})$,
\item"(B)" $(\Bbb{P}^{n-1}, T_{\Bbb{P}^{n-1}})$,
\item"(C)" $(\Bbb{Q}^{n-1}, \Cal{O}_{\Bbb{Q}^{n-1}}(1)^{\oplus n-1})$,
\item"(D)" $F$ is a $\Bbb{P}^{n-2}$-bundle over a smooth curve $B$
           and $\Cal{E}_{F_{\pi}}=\Cal{O}_{\Bbb{P}^{1}}(1)^{\oplus n-1}$
           for every fiber $F_{\pi}$ of the projection $\pi: F\to B$.
\endroster
If $(F,\ce_{F})$ is one of the type (A), (B), or (C),
then $h^{0}(K_{F}+c_{1}(\ce_{F}))>0$ by easy calculation.
Here we prove the following Lemma.

\proclaim{Lemma 5.5} 
Let $(f,X,C,\ce)$ be a generalized polarized fiber space
with $n\geq r+1$ and $\dim C=1$,
where $n:=\dim X\geq 3$ and $r:=\rk\ce\geq 2$.
If $h^{0}(K_{F}+c_{1}(\ce_{F}))>0$ 
for a general fiber $F$ of $f$, 
then $(K_{X/C}+c_{1}(\ce))c_{1}(\ce)^{n-r-1}c_{r}(\ce)>0$.
\endproclaim

\demo{Proof}
By hypothesis, $f_{*}\Cal{O}(K_{X/C}+c_{1}(\ce))\not=0$.
By Theorem 2.4 and Corollary 2.5 in \cite{EV},
we get that $f_{*}\Cal{O}(K_{X/C}+c_{1}(\ce))$ is ample.
By the proof of Lemma 1.4.1 in \cite{Fk2},
we get that $m(K_{X/C}+c_{1}(\ce))-f^{*}A$ is an effective divisor
for a large number $m$ and an ample divisor $A$ on $C$.
Hence we obtain $(K_{X/C}+c_{1}(\ce))c_{1}(\ce)^{n-r-1}c_{r}(\ce)>0$. 
\qed
\enddemo
 
By Lemma 5.5, we get that $g(X,\ce)>g(C)$ if 
$(F,\ce_{F})$ is one of the type (A), (B), or (C).

Assume that $(F, \ce_{F})$ is the type (D).
Then there exist vector bundles $\cf$ and $\cg$ on $B$ with 
$\rk\cf=\rk\cg=n-1$
such that $\ce_{F}\cong H(\cf)\otimes\pi^{*}(\cg)$,
where $H(\cf)$ is the tautological line bundle of $\pp(\cf)$.
Then $K_{F}+c_{1}(\ce_{F})=\pi^{*}(K_{B}+\det\cf+\det\cg)$.
Since $K_{F}+c_{1}(\ce_{F})$ is nef, we get 
$(K_{X/C}+c_{1}(\ce))c_{r}(\ce)\geq 0$ by the proof of Lemma 5.4.
We have $g(X,\ce)=g(C)$, then $c_{r}(\ce)F=1$.
Since $1=c_{r}(\ce_{F})=c_{1}(\cf)+c_{1}(\cg)$,
we obtain that
$$\split
&h^{0}(K_{B}+\det\cf+\det\cg) \\
&\geq 1-g(B)+\deg(K_{B}+\det\cf+\det\cg) \\
&=g(B)-1+c_{1}(\cf)+c_{1}(\cg) \\
&=g(B).
\endsplit
$$
Because $K_{F}+c_{1}(\ce_{F})$ is nef, we obtain that
$\deg(K_{B}+\det\cf+\det\cg)\geq 0$.
Hence $g(B)\geq 1$.
Therefore $h^{0}(K_{F}+c_{1}(\ce_{F}))\geq 1$.
By Lemma 5.5 we obtain that $g(X,\ce)>g(C)$ and this is a contradiction.

Next we assume that $K_{F}+c_{1}(\ce_{F})$ is not nef.
Then $K_{X}+c_{1}(\ce)$ is not nef, and by \cite{YZ},
there exists a surjective morphism $\pi : X\to B$ such that any fiber 
of $\pi$ is $\pp^{n-1}$ and $\ce|_{F_{\pi}}\cong \oplus^{n-1}\Cal{O}(1)$,
where $B$ is a smooth curve and $F_{\pi}$ is a fiber of $\pi$.
\flushpar
Since any fiber of $\pi$ is $\pp^{n-1}$, there exists a morphism
$\delta: B\to C$ such that $f=\delta\circ\pi$.
Because $f$ has connected fibers, $\delta$ is an isomorphism.
This completes the proof of Theorem 5.3. \qed
\enddemo

\example{Remark 5.6}
  Let $(f,X,C,\ce)$ be as in Theorem 5.2.
  Suppose that $g(X,\ce)=g(C)$ and $r=1$.
  Then by Theorem 1.4.2 and Proposition 1.4.3 in \cite{Fk2},
  $(f,X,C,\ce)$ is a scroll unless $n=2$ and 
  $(f,X,C,\ce)\cong(\pi,\pp^1\times\pp^1,\pp^1,L)$,
  where $\pi$ is one projection such that $LF_{\pi}\ge 2$
  for a fiber $F_{\pi}$ of $\pi$.
\endexample

\subhead{Appendix}
\endsubhead

\proclaim{Proposition A}
Let $(X,L)$ be a quasi-polarized surface 
(i.e. $L$ is a nef and big line bundle on a smooth surface $X$)
such that $\kappa(X)=2$ and $h^{0}(L)\geq 2$.
Then $K_{X}L\geq 2q(X)-2$.
If equality holds and $(X,L)$ is $L$-minimal
(i.e. $LE>0$ for any $(-1)$-curve $E$ on $X$),
then $(X,L)$ is the following:
\flushpar
$X\cong F\times C$ and $L\equiv C+2F$,
where $F$ and $C$ are smooth curves with $g(F)=2$ and $g(C)\geq 2$.
\endproclaim

\demo{Proof}
See \cite{Fk4}.\qed
\enddemo

\proclaim{Proposition B}
Let $(X,L)$ be a polarized surface with $\kappa(X)=0$ or $1$.
Assume that $L$ is spanned.
Then $g(L):=g(X,L)\geq 2q(X)$.
Furthermore if $g(L)=2q(X)$, then $(X,L)$ is one of the following:
\roster
\item $(X,L)$ is a polarized abelian surface with $L^{2}=6$
      such that $(X,L)\not\cong 
      (E_{1}\times E_{2}, p_{1}^{*}(D_{1})+p_{2}^{*}(D_{2}))$,
      where $E_{i}$ is a smooth elliptic curve,
      $p_{i}$ is the $i$-th projection,
      and $D_{i}\in\operatorname{Pic}(E_{i})$ for $i=1,2$ 
      with $\deg D_{1}=1$ and $\deg D_{2}=3$.
\item $X$ is a one point blowing up of $S$, and $L=\mu^{*}A-2E$,
      where $S$ is an abelian surface, $A$ is an ample line bundle with 
      $A^{2}=8$, $\mu: X\to S$ is its blowing up, 
      and $E$ is a $(-1)$-curve of $\mu$.
\item $\kappa(X)=1$, $L^{2}=4$, $q(X)=3$, $X$ has a locally trivial 
      elliptic fibration
      $f: X\to C$, and $LF=3$ for a fiber $F$ of $f$,
      where $C$ is a smooth curve with $g(C)=2$.
\endroster
\endproclaim

\demo{Proof}
See \cite{Fk5}. \qed
\enddemo

\proclaim{Theorem}
Let $X$ be a smooth projective surface and let $L$ be an ample and 
spanned line bundle on $X$.
If $g(L)=q(X)+2$, then $(X,L)$ is one of the following:
\roster
\item $(X,L)$ is a relatively minimal conic bundle over a smooth curve
      $B$ of genus two (i.e. $X$ is a $\pp^1$-bundle over $B$ and
      $L_F=\co_{\pp^1}(2)$ for every fiber $F$ of the ruling).
\item $X$ is a $\Bbb{P}^{1}$-bundle $X_{0}$
      blown-up at $s$ $(0\le s\le 4)$ points $p_{1}, \dots, p_{s}$ 
      on distinct fibers and $L=\pi^{*}L_{0}-E_{1}-\dots -E_{s}$, 
      where $\pi: X\to X_{0}$ is the blowing up, $E_i=\pi^{-1}(p_i)$,
      $X_{0}$ is an elliptic $\Bbb{P}^{1}$-bundle of invariant $e\leq 0$,
      and $L_{0}\equiv 2\sigma+(e+2)f$ 
      ($\sigma$ is a minimal section with $\sigma^2=-e$
      and $f$ is a fiber). 
\item $X$ is an $\Bbb{F}_{e}$ $(e\leq 2)$ blown-up at $s$ $(0\le s\le 9)$
      points $p_{1}, \dots, p_{s}$ on distinct fibers and
      $L=\pi^{*}L_{0}-E_{1}-\dots -E_{s}$,
      where $\pi: X\to \Bbb{F}_{e}$ is the blowing up, 
      $E_{i}=\pi^{-1}(p_{i})$, and $L_{0}\equiv 2\sigma+(e+3)f$.
\item $X$ is a del Pezzo surface of degree one and there exists 
      a double covering $\pi: X\to\Cal{Q}\subset\Bbb{P}^{3}$
      of a quadric cone $\Cal{Q}$
      branched at the vertex and along the transverse intersection of 
      $\Cal{Q}$ with a cubic surface and $L=\pi^{*}(\Cal{O}_{\Cal{Q}}(1))$.
\item $(X,L)$ is a polarized abelian surface with $L^{2}=6$
      such that $(X,L)\not\cong 
      (E_{1}\times E_{2}, p_{1}^{*}(D_{1})+p_{2}^{*}(D_{2}))$,
      where $E_{i}$ is a smooth elliptic curve,
      $p_{i}$ is the $i$-th projection,
      and $D_{i}\in\operatorname{Pic}(E_{i})$ for $i=1,2$ 
      with $\deg D_{1}=1$ and $\deg D_{2}=3$.
\item $X$ is a blowing up of an abelian surface $S$ at one point $p$
      and $L=\pi^{*}A-2E$,
      where $\pi:X\to S$ is the blowing up, $E=\pi^{-1}(p)$, 
      and $A$ is an ample line bundle on $S$ with $A^{2}=8$. 
\item $X$ is a $K3$ surface which is a double covering of $\Bbb{P}^{2}$ 
      branched along a smooth curve of degree six and 
      $L$ is the pull back of $\Cal{O}_{\Bbb{P}^{2}}(1)$.
\endroster
\endproclaim

\demo{Proof}
(I) The case in which $\kappa(X)=0$ or 1.
\flushpar
Then by Proposition B, we get that
$g(L)\geq 2q(X)$.
So we obtain $q(X)\leq 2$ by assumption.
\flushpar
(I-1) If $q(X)=2$, then $g(L)=q(X)+2=2q(X)$ and by Proposition B
we get the type (5) and (6) in Theorem.
\flushpar
(I-2) If $q(X)\leq 1$, then $g(L)\leq 3$ and $L^{2}\leq 4$ by $K_XL\ge 0$
\flushpar
(I-2-1) If $L^{2}=4$, then $\kappa(X)=0$ and $X$ is minimal 
since $K_{X}L=0$.
So by Kodaira vanishing Theorem and Riemann-Roch Theorem,
we get the equality:
$h^{0}(L)=L^{2}/2+\chi(\Cal{O}_{X})=2+\chi(\Cal{O}_{X})$.
Because $L$ is ample and spanned, we obtain $h^{0}(L)\geq 3$
and $\chi(\Cal{O}_{X})\geq 1$.
But then $q(X)=0$ by the classification theory of surfaces
and this is impossible.
\flushpar
(I-2-2) If $L^{2}=3$, then $g(L)=3$, $K_{X}L=1$, and $q(X)=1$.
We have $h^0(L)\ge 3$ since $L$ is ample spanned.
\flushpar
If $h^{0}(L)\geq 4$, then $g(L)>\Delta(L)$ and $L^{2}\geq 2\Delta(L)+1$,
where $\Delta(L):=2+L^{2}-h^{0}(L)$ is the $\Delta$-genus of $L$.
But then $q(X)=0$ (see e.g. (I.3.5) in \cite{Fj4}).
\flushpar
If $h^{0}(L)=3$, then there is a triple covering 
$\varphi_{|L|}: X\to \Bbb{P}^{2}$ which is defined by $|L|$.
Let $\Cal{E}$ be a vector bundle of rank two on $\Bbb{P}^{2}$
such that $\pi_{*}\Cal{O}_{X}=\Cal{O}_{\Bbb{P}^{2}}\oplus\Cal{E}$.
By Lemma 3.2 in \cite{Be}, we get the following two equalities:
\roster
\item"(i)" $\chi(\Cal{O}_{X})=(1/2)g(L)(g(L)+1)+2-c_{2}$,
\item"(ii)" $K_{X}^{2}=2g(L)^{2}-4g(L)+11-3c_{2}$,
\endroster
where $c_{2}:=c_{2}(\Cal{E})$.
Since $g(L)=3$, we get that $3\chi(\Cal{O}_{X})-K_{X}^{2}=7$
by the above equalities.
\flushpar
If $\kappa(X)=0$, then $K_{X}^{2}=-1$ because $K_{X}L=1$.
So we get $\chi(\Cal{O}_{X})=2$.
But by the classification theory of surfaces, 
this is impossible because $q(X)=1$.
\flushpar
If $\kappa(X)=1$, then $X$ is minimal and $K_{X}^{2}=0$
because $K_{X}L=1$.
But then $3\chi(\Cal{O}_{X})=7$ and this is impossible.
\flushpar
(I-2-3) If $L^{2}=2$, then $K_{X}L=0$ or 2.
Since $\kappa(X)\geq 0$, we get that $\Delta(L)\ge 1$ and $h^{0}(L)=3$.
Then there exists a double covering 
$\varphi_{|L|}: X\to \Bbb{P}^{2}$ which is defined by $|L|$.
We remark that $K_{X}=\varphi_{|L|}^{*}(K_{\Bbb{P}^{2}}+D)$
for some $D\in \operatorname{Pic}(\Bbb{P}^{2})$.
Since $\kappa(X)=0$ or 1, we get that $\kappa(X)=0$ and so $X$ is minimal.
In particular $K_{X}=\Cal{O}_{X}$.
Therefore $K_{X}L=0$ and $g(L)=2$.
Since $h^{0}(L)=L^{2}/2+\chi(\Cal{O}_{X})=1+\chi(\Cal{O}_{X})$,
we get $\chi(\Cal{O}_{X})=2$.
Hence $X$ is a K3 surface by the Classification theory of surfaces.
This is the type (7) in Theorem.
\flushpar
(II) The case in which $\kappa(X)=2$.
\flushpar
Then by Corollary 3.2 in \cite{Fk1}, we get
$g(L)\geq 2q(X)-1$.
So we obtain $q(X)\leq 3$ and $g(L)\leq 5$ by assumption.
Furthermore $L^{2}\leq 3$ by Proposition A
because $L$ is spanned.
(We remark that $L$ is $L$-minimal if $L$ is ample.)
\flushpar
If $h^{0}(L)\geq 4$, then $g(L)>1\geq \Delta(L)$ and
$L^{2}\geq 2\Delta(L)+1$.
On the other hand, since $\kappa(X)\geq 0$, we obtain that
$\Delta(L)=1$ and $L^{2}=3$.
So we get $q(X)=0$ and $g(L)\geq 3$ and this is impossible.
Therefore $h^{0}(L)=3$.
\flushpar
If $L^{2}=3$, then there exists a 
triple covering 
$\varphi_{|L|}: X\to \Bbb{P}^{2}$ which is defined by $|L|$.
In this case, by the same argument as above,
we get 
$$2(K_{X}^{2}-3\chi(\Cal{O}_{X}))=(g(L)-1)(g(L)-10).$$
Since $3\leq g(L)\leq 5$, we get the following:

\roster
\item"($\alpha$)"
 $(g(L), q(X), K_{X}L, K_{X}^{2}-3\chi(\Cal{O}_{X}))
=(5,3,5,-10)$,
\item "($\beta$)"
$(g(L), q(X), K_{X}L, K_{X}^{2}-3\chi(\Cal{O}_{X}))
=(4,2,3,-9)$,
\item "($\gamma$)"
$(g(L), q(X), K_{X}L, K_{X}^{2}-3\chi(\Cal{O}_{X}))
=(3,1,1,-7)$.
\endroster

\proclaim{Claim}
The above three cases cannot occur.
\endproclaim

\demo{Proof}
(II-1) The case $(\gamma)$.
\flushpar
In this case $X$ is minimal because $K_{X}L=1$.
But then this is impossible by Hodge index Theorem.
\flushpar
(II-2) The case $(\beta)$.
\flushpar
If $X$ is minimal, then $K_{X}^{2}\geq 2q(X)=4$ by Th\'eor\`eme 6.1
in \cite{De}.
On the other hand, $K_{X}^{2}\leq 3$ by Hodge index Theorem
and this is a contradiction.
\flushpar
So we get that $X$ is not minimal.
Let $\mu:=\mu_{r}\circ\dots\circ\mu_{1}: X:=X_{0}\to X_{1}\to\dots
\to X_{r-1}\to X_{r}=:X^{\prime}$ be an admissible minimalization of $X$ 
and let $m=(m_{r},\dots , m_{1})$ be the weight sequence of this 
minimalization (see (II.14.4) in \cite{Fj4}).
We remark that $m_{r}\geq \dots \geq m_{1}$.
\flushpar
If $m_{1}=1$, then $g(L_{1})=q(X_{1})+1$ and $h^{0}(L_{1})\geq 2$,
where $L_{1}:=(\mu_{1})_{*}(L)$ in the sense of cycle theory.
But then this is impossible by Proposition A 
because $2=K_{X}L>K_{X_{1}}L_{1}$.
So we get $m_{1}\geq 2$.
Then $L_{1}^{2}\geq 7$ and $K_{X_{1}}L_{1}\leq 1$.
Hence $X_{1}$ is minimal and this is a contradiction by Hodge index Theorem.
\flushpar
(II-3) The case $(\alpha)$.
\flushpar
If $X$ is minimal, then $\chi(\Cal{O}_{X})\geq 4$ because
$3\chi(\Cal{O}_{X})=K_{X}^{2}+10$.
Furthermore $p_{g}(X)\geq 6$ since $q(X)=3$.
Hence $K_{X}^{2}\geq 2p_{g}(X)\geq 12$ by Th\'eor\`eme 6.1 in \cite{De}.
But this is impossible by Hodge index Theorem.
So we get that $X$ is not minimal.
By the same argument as in the case (II-2) we get a contradiction.\qed
\enddemo
\flushpar
We continue the proof of Theorem.
\flushpar
If $L^{2}=2$, then there exists a double covering 
$\varphi_{|L|}: X\to \Bbb{P}^{2}$ which is defined by $|L|$.
Let $\Cal{O}_{\Bbb{P}^{2}}(a)$ be 
a line bundle on $\Bbb{P}^{2}$ such that 
$B\in |\Cal{O}_{\Bbb{P}^{2}}(2a)|$,
where $B$ is the branch locus.
Then $(\varphi_{|L|})_{*}(\Cal{O}_{X})=\Cal{O}_{\Bbb{P}^{2}}\oplus
\Cal{O}_{\Bbb{P}^{2}}(-a)$.
Hence
$$
h^{1}(\Cal{O}_{X})=h^{1}((\varphi_{|L|})_{*}(\Cal{O}_{X})) 
                  =h^{1}(\Cal{O}_{\Bbb{P}^{2}})
                   +h^{1}(\Cal{O}_{\Bbb{P}^{2}}(-a)) 
                  =0.
$$
So we get $g(L)=2$.
But since $K_{X}L>0$ and $L^{2}=2$, this is impossible.
\flushpar
(III) The case in which $\kappa(X)=-\infty$.
\flushpar
Since $(X,L)$ is not a scroll over a smooth curve,
we get $g(L)\geq 2q(X)$ by Claim 3.2.
So $q(X)\leq 2$. 
\flushpar
(III-1) The case in which $q(X)=2$.
\flushpar
In this case, $g(L)=q(X)+2=2q(X)$.
As in Claim 3.2, we obtain that 
$$0\leq (K_{X}+L)^{2}\leq 4(g(L)-2q(X)+1)-L^{2}.$$
Hence $L^{2}\leq 4$ in this case.
\flushpar
If $L^{2}=4$, then $X$ is relatively minimal and $(K_{X}+L)^{2}=0$,
that is, $(X,L)$ is a relatively minimal conic bundle over a smooth curve.
This is the type (1) in Theorem.
\flushpar
If $L^{2}\leq 3$ and $h^{0}(L)\geq 4$, then we get a contradiction
as in (I-2-2).
So we may assume that $L^{2}\leq 3$ and $h^{0}(L)=3$.
\flushpar
If $L^{2}=3$, then $K_{X}L=3$ and
there is a triple covering 
$\varphi_{|L|}: X\to \Bbb{P}^{2}$ which is defined by $|L|$.
Since $\chi(\Cal{O}_{X})=-1$, we get that $K_{X}^{2}=-12$
by Lemma 3.2 in \cite{Be}.
Here we calculate $(K_{X}+L)^{2}$;
$$ (K_{X}+L)^{2}=K_{X}^{2}+2K_{X}L+L^{2}=-12+6+3<0. $$
But this is a contradiction because $K_{X}+L$ is nef.
\flushpar
If $L^{2}=2$, then 
there is a double covering 
$\varphi_{|L|}: X\to \Bbb{P}^{2}$ which is defined by $|L|$.
But then $q(X)=0$ and this is a contradiction.
\flushpar
(III-2) The case in which $q(X)=1$.
\flushpar
Then $g(L)=3$.
Here we use the classification of polarized surfaces
with sectional genus three by \cite{LL}.

\proclaim{Claim}
The case in which $L^{2}=3$ cannot occur.
\endproclaim

\demo{Proof} If $L^{2}=3$ and $h^{0}(L)\geq 4$, then 
$g(L)>1\geq \Delta(L)$ and $L^{2}\geq 2\Delta(L)+1$.
But this is impossible because $q(X)=1$.
So we may assume that $h^{0}(L)=3$.
Then there is a triple covering 
$\varphi_{|L|}: X\to \Bbb{P}^{2}$ which is defined by $|L|$.
Since $\chi(\Cal{O}_{X})=0$, we get $K_{X}^{2}=-7$ by 
Lemma 3.2 in \cite{Be}.
But in the table II of \cite{LL},
the case in which $L^{2}=3$ cannot occur. \qed
\enddemo

\flushpar
Next we prove that the following case cannot occur
(see (2.6) in \cite{LL}):
\flushpar
{\it $X$ is an elliptic $\Bbb{P}^{1}$-bundle $X_{\sharp}$ 
of invariant $e=0$, blown up at a single point $p$ not lying on 
a curve $D\in |m\sigma|$, $m\leq 2$ and 
$L=\eta^{*}[4\sigma+(2e+1)f]\otimes [E]^{-2}$.}
(Here we use the same notations as in \cite{LL}.)
\flushpar
Let $\sigma'$ be the strict transform of $\sigma$ under $\eta$.
Since 
$$0<L\sigma'=(4\sigma+f)\sigma-2E\sigma'=1-2E\sigma',$$
we see that $E\sigma'=0$ and $L\sigma'=1$.
It follows that $\sigma\cong\sigma'\cong\Bbb{P}^{1}$
since $L$ is spanned.
This is a contradiction.
\flushpar
By the above argument, we obtain the type (2) in Theorem
by the classification of polarized surfaces with sectional genus three
(see \cite{LL}).
\flushpar
(III-3) The case in which $q(X)=0$.
\flushpar
Then $g(L)=2$.
So by Theorem 3.1 in \cite{LP} we get the type (3) and (4)
in Theorem. \qed
\enddemo

\enddocument